\documentclass{crmp-l}

\usepackage{amsmath}
\usepackage{amsfonts}
\usepackage{amssymb}
\usepackage{amstext}
\usepackage{amsbsy}
\usepackage{amsthm}
\usepackage{amscd}
\usepackage[all]{xy}

\newtheorem{theorem}{Theorem}[section]

\newtheorem{corollary}[theorem]{Corollary}

\theoremstyle{definition}

\theoremstyle{remark}

\numberwithin{equation}{section}

\newcommand{\bb}[1]{\ensuremath{\mathbb{#1}}}
\newcommand{\dd}[3]{\ensuremath{\frac{\mathrm{d}^{#3} {#1}}{\mathrm{d} {#2}^{#3}}}}

\title{Normal Forms, K3 Surface Moduli, and Modular Parametrizations}

\author{A.~Clingher}
\address{Department of Mathematics and Computer Science, University of Missouri -- St.~Louis, St.~Louis, MO 63121}
\email{clingher@arch.umsl.edu}
\thanks{A.C.~acknowledges support from a University of Missouri Research Board Grant.  He would also like to thank
the Department of Mathematics of the University of Washington for its hospitality and financial support during
the 2007 summer school on K3 Surfaces and String Duality.}

\author{C.~F.~Doran}
\address{Department of Mathematics, University of Washington, Seattle, WA 98195}
\email{doran@math.washington.edu}
\thanks{C.F.D. was supported in part by the University of Washington Royalty Research Fund.}

\author{J.~Lewis}
\address{Department of Mathematics, University of Washington, Seattle, WA 98195}
\email{jacobml@math.washington.edu}
\thanks{J.L.~was supported in part by a VIGRE Graduate Fellowship.}

\author{U.~Whitcher}
\address{Department of Mathematics, University of Washington, Seattle, WA 98195}
\email{ursula@math.washington.edu}
\thanks{U.W.~was supported in part by a McKibben and Merner Fellowship.}

\dedicatory{For John McKay}

\begin{document}

\maketitle

\section{Introduction}

The geometric objects of study in this paper are K3 surfaces which admit a polarization by the unique even
unimodular lattice of signature (1,17). This lattice can be realized as the orthogonal direct sum
$${\rm M} = {\rm H} \oplus
{\rm E}_8 \oplus {\rm E}_8 \ ,$$ where ${\rm H}$ is the standard rank-two hyperbolic lattice and ${\rm E}_8$ is
the unique even, negative-definite and unimodular lattice of rank eight. A lattice-polarization by the lattice
${\rm M}$ imposes severe constraints on the geometry of a K3 surface $X$.  In particular, the Picard rank of $X$
has to be 18, 19, or 20 (see \cite{Mor}).

A standard Hodge-theoretic observation about this special class of K3 surfaces is that the polarized Hodge
structure of an ${\rm M}$-polarized K3 surface is identical with the polarized Hodge structure of an abelian
surface $A = {\rm E}_1 \times {\rm E}_2$ realized as a cartesian product of two elliptic curves.  Since both
types of surfaces involved admit appropriate versions of the Torelli theorem, Hodge theory implies a
well-defined correspondence, giving rise to a canonical analytic isomorphism between the moduli spaces on the
two sides.  By employing a modern point of view from the frontier of algebraic geometry with string theory, one
can regard this correspondence as a Hodge-theoretic {\bf duality map}, a correspondence that relates two
seemingly different types of surfaces sharing similar Hodge-theoretic information.

In \cite{CD1}, the first two authors showed that that the resemblance of the two Hodge structures involved in
the duality correspondence is not fortuitous, but rather is merely a consequence of a quite interesting
geometric relationship.   We show in \cite[Theorem 1.1]{CD1} that: the surface $X$ possesses a canonical
involution defining a Shioda-Inose structure; the minimal resolution of $X$ quotiented by this involution is a
new K3 surface $Y$ endowed with a canonical Kummer structure; this structure realizes $Y$ as the Kummer surface
of an abelian surface $A$ which is canonically represented as a Cartesian product of two elliptic curves (i.e.,
$A = {\rm E}_1 \times {\rm E}_2$); and the construction induces a canonical Hodge isomorphism between the ${\rm
M}$-polarized Hodge structure of $X$ and the natural ${\rm H}$-polarized Hodge structure of the abelian surface
$A$.

One important feature of the special class of ${\rm M}$-polarized K3 surfaces is that they turn out to be
completely classified by two modular invariants $\sigma$ and $\pi$ in ${\mathbb C}$, much in the same way as
elliptic curves over the field of complex numbers are classified by the $j$-invariant.  However, the two modular
invariants $\sigma$ and $\pi$ are not geometric in origin.  They are defined Hodge-theoretically, and the result
leading to the classification is a consequence of the appropriate version of the Global Torelli Theorem for
lattice polarized K3 surfaces. In the context of the duality map, the two invariants can be seen as the standard
symmetric functions on the $j$-invariants of the dual elliptic curves, $\sigma$ being their sum and $\pi$ their
product.  This interpretation suggests that the modular invariants of an ${\rm M}$-polarized K3 surface can be
computed by determining the two elliptic curves.

Explicit ${\rm M}$-polarized K3 surfaces can be constructed by various geometrical procedures.  One such method,
introduced in 1977 by Inose \cite{inose1}, constructs a two-parameter family $X(a,b)$ of ${\rm M}$-polarized K3
surfaces by taking minimal resolutions of the projective quartics in ${\mathbb P}^3$ associated with the special
equations:
\begin{equation}
\label{InoseIntro}
 y^2 z w - 4 x^3 z + 3 a x z w^2 - \frac{1}{2} (z^2 w^2 + w^4) + b z w^3 = 0
\end{equation}
for $a, b$ in ${\mathbb C}$. In Section \ref{K3PolHE8E8} below, this construction is generalized so as to
compactify the affine $(a,b)$-parametric family to a weighted projective space with covering coordinates $[a, b,
d]$.  The coarse moduli space for ${\rm M}$-polarized K3 surfaces is then given by the locus where $d \neq 0$ in
this weighted projective space.  One can thus regard the Inose quartic as a {\bf normal form for ${\rm
M}$-polarized K3 surfaces}. In the second part of \cite{CD1}, we use the geometric correspondence to explicitly
describe the $j$-invariants of the two elliptic curves ${\rm E}_1$ and ${\rm E}_2$ associated to the Inose
surface $X(a,b)$. The result is reproduced below (Theorem \ref{theo3}) in the context of the $X(a,b,d)$ family.
As a consequence, the coarse moduli space for ${\rm M}$-polarized K3 surfaces is explicitly identified with the
classical Hilbert modular surface $\mathcal{H}_1$ for the elliptic modular group.

In this paper, we use explicit computations of Picard-Fuchs differential equations to explore (and exploit) the
``universal'' property of the moduli space of ${\rm M}$-polarized K3 surfaces with respect to enhancement to
{\em non-generic} Picard lattices.  The conditions on the hypersurfaces for enhancement to Picard lattices of
rank 19 are of particular interest, as these conditions determine one-parameter subfamilies of the coarse moduli
space.  As we work throughout in parallel with the case of the Weierstrass normal form for elliptic curves, this
case is reviewed in Section \ref{WeirEC}.  The approach is a natural extension of the second author's work on
Picard-Fuchs uniformization for families of elliptic curves and lattice-polarized K3 surfaces \cite{Dor1,Dor2}.

The Inose family normal form is realized as projective hypersurfaces with just ADE type singularities.  This
provides a crucial advantage to using this normal form, as now the Griffiths-Dwork technique (see Section
\ref{GDGeneral} for a general discussion) can be applied to compute the differential equations satisfied by
periods of the holomorphic 2-forms directly. We do this first in Section 3 with independent variables taken to
be the free coefficients in the normal form hypersurface equations; this is the analog of taking partial
derivatives of the Weierstrass elliptic periods with respect to the parameters $g_2$ and $g_3$.  It proves
convenient to work in the affine chart $[1,b,d]$ for these computations.  Upon substitution of the modular
invariants from \cite{CD1}, this partial differential system completely decouples, reducing to a pair of ODEs.
Each of these is the Picard-Fuchs ODE of a corresponding ``factor'' family of elliptic curves, reflecting the
geometric fact that our normal form K3 surfaces are Shioda-Inose surfaces of products of pairs of elliptic
curves.

Next we apply the Griffiths-Dwork method with independent variable $t$; this is the analog of treating $g_2(t)$
and $g_3(t)$ as rational functions of a single complex variable $t$. In the case of Weierstrass normal form
elliptic curves, this modification simply allows one to compute the Picard-Fuchs ODE for arbitrary one parameter
families of elliptic curves, as the moduli space of elliptic curves is itself one-dimensional.  In our case,
this modification allows one to ``probe'' the two-dimensional moduli space of ${\rm H}\oplus {\rm E}_8\oplus
{\rm E}_8$ lattice polarized K3 surfaces with rational curves, computing the Picard-Fuchs ODE for the restricted
period functions.\footnote{Computation of Picard-Fuchs equations for particular one-parameter families of K3
surfaces of Picard rank 19 dates back at least to \cite{Pet}.  The approach taken in the present paper, however,
results in a universal expression for the Picard-Fuchs equations for {\em all} one-parameter families of K3
surfaces admitting (at least) an $M$-polarization.  The families are explicitly realized as subfamilies of the
single two-parameter family of Equation~\ref{InoseIntro}.} The modular invariants of \cite{CD1} can once again
be substituted in, yielding a fourth order ODE whose coefficients are differential rational functions of
$j_1(t)$ and $j_2(t)$.

Here is where the universal property for the moduli space with respect to lattice enhancement makes its impact
felt.  Enhancement of polarization from ${\rm M}$ to
$${\rm M}_n = {\rm H}\oplus {\rm E}_8\oplus {\rm E}_8\oplus \langle - 2 n \rangle$$ is
equivalent to the condition that this fourth order ODE reduces to a third order ODE (a condition on ODEs studied
by Gino Fano over a century ago). Under the hypothesis that $j_1(t)$ and $j_2(t)$ are nonconstant, this
reduction occurs precisely when these two functions satisfy a particular nonlinear ODE.  This nonlinear ODE is
the {\bf master equation for modular parametrizations} in that any parametrization of the modular curve $X_0(n)$
by rational functions $j_1(t)$ and $j_2(t)$ solves this ODE.  What's more, for any solutions $j_1(t)$ and
$j_2(t)$ not themselves rational functions of $t$, but such that their sum and product are still rational
functions, the pair parametrizes the modular curve $X_0(n)+n$.

The master equation can even be rewritten in a form which makes manifest its relationship with pairs of
$n$-isogenous elliptic curves.  First, all of the terms involving $j_1(t)$ can be moved to one side of the
equation and all those involving $j_2(t)$ to the other.  Of course, up to exchanging $j_1(t)$ with $j_2(t)$, the
expressions on each side are then the same.  Moreover, they are recognizable as the coefficients of projective
normalized Picard-Fuchs ODEs for a pair of elliptic curves. As such, they are each the sum of one term involving
the rational function characterizing the uniformizing differential equation for the elliptic modular parameter
$j$, and another term consisting of the Schwarzian derivative of $j$ with respect to $t$ (when set equal to
zero, this expression is also known as the ``Schwarzian differential equation'' for $j$). Once the master
equation for modular parametrizations has been expressed in this form, it is natural to ask whether an analogous
equation, based on hauptmoduls other than the elliptic modular parameter $j$, will similarly characterize
parametrizations of genus zero quotients of the upper half plane by the moonshine groups. We address this
question in Section \ref{hauptmoduls}, and illustrate there explicit parametrizations of some modular equations
studied by Cohn-McKay.

Currently under investigation are an extension of the method to handle modular curves of genus greater than zero
(by modifying the Griffiths-Dwork technique) and a generalization to parametrizing Humbert surfaces and Shimura
curves based now on universal properties of the moduli space of ${\rm H}\oplus {\rm E}_8\oplus {\rm E}_7$
lattice polarized K3 surfaces studied by Clingher and Doran in \cite{CD2,CD3}.

\section{Elliptic Curves}
\label{WeirEC}

\subsection{Normal Form and Moduli Space}

As is well-known, the classification of elliptic curves over the field of
complex numbers is based on the following
classical features. First, every elliptic curve ${\rm E}$ can be realized
through an explicit {\bf normal form} given,
for instance, by the projective version of the Weierstrass cubic:
\[ y^2 z - 4 x^3 + g_2 x z^2 + g_3 z^3 = 0\]
for some pair $(g_2,g_3) \in \mathbb{C}^2$ with $g_2^3-27g_3^2 \neq 0$. The Weierstrass form is not unique.
Nevertheless, the weighted projective point $[g_2,g_3] \in \mathbb{WP}(2,3)$ classifies ${\rm E}$ up to
isomorphism. This leads to a nice (coarse) {\bf moduli space} for complex elliptic curves, taken as the
one-dimensional open variety:
$$ \mathcal{M} \ = \ \left \{ \
[g_2,g_3] \in \mathbb{WP}(2,3) \ \middle \vert \ g_2^3-27g_3^2 \neq 0
\ \right \}. $$
In this context, the {\bf {\it j}-invariant}:
\[ j({\rm E}) \ = \ \frac{g_2^3}{g_2^3-27g_3^2} \]
realizes an identification between $\mathcal{M}$ and $\mathbb{C}$.
\par From a Hodge theoretic point of view, an elliptic curve can be seen
as an element (period) in the
complex upper half-plane $\mathbb{H}$. The classifying space for such
periods is the {\bf modular quotient}
$$ \mathcal{F}_1 \ = \ \Gamma_1 \backslash \mathbb{H} $$
where $\Gamma_1 = {\rm PSL}(2, \mathbb{Z})$ with the standard action on
$\mathbb{H}$.
\par The link between the two points of view is realized by the period map:
\begin{equation}
\label{inverseper1}
 {\rm per} \colon \mathcal{M} \ \rightarrow \ \mathcal{F}_1
\end{equation}
which makes an isomorphism of analytic spaces. The inverse  of
$(\ref{inverseper1})$ can be described in terms
of $\Gamma_1$-modular forms as:
\[ {\rm per}^{-1} \ = \ \left [ \ 60 \cdot {\rm E}_4, \ 140 \cdot {\rm E}_6
\ \right ] \] where ${\rm E}_4, {\rm E}_6 \colon \mathbb{H} \rightarrow \mathbb{C}$ are the classical Eisenstein
series of weights four and six, respectively.  Using the inverse period map, one may view the $j$-invariant as a
map $j: \mathbb{H} \to \mathbb{C}$.

\subsection{Picard--Fuchs Differential Equation}

When considering the periods of varieties in families, it is often easier to compute the regular singular
differential equation satisfied by all the period functions on the base of the family -- the {\bf Picard-Fuchs
differential equation} -- than to explicitly describe the period functions themselves.  We first review in
Section \ref{GDGeneral} the method, due to Griffiths and Dwork, for computing the Picard-Fuchs differential
equation for families of hypersurfaces in projective space via residues.  In Section \ref{GDforCurves} this
technique is applied to families of elliptic curves in ${\mathbb P}^2$ in Weierstrass normal form.

\subsubsection{The Griffiths--Dwork Technique} \label{GDGeneral}

Let $X$ be a hypersurface in $\mathbb{P}^n$ given by a homogeneous polynomial $Q$ in coordinates $[x_0, \dots,
x_n]$, and let $\iota: X \to \mathbb{P}^n$ be the inclusion map. Let $\mathcal{H}(X)$ be the de Rham cohomology
of rational $n$-forms on $\mathbb{P}^n - X$.  We may write any representative of $\mathcal{H}(X)$ as $P
\Omega_0/Q^k$, where $\Omega_0=\sum_{i=0}^n (-1)^i x^i \mathrm{d}{x^0}\wedge \ldots
\widehat{\mathrm{d}{x^i}}\ldots\wedge\mathrm{d}{x^n}$ is the usual holomorphic form on $\mathbb{P}^n$ and $P$ is
a homogeneous polynomial of degree $\mathrm{deg} \,P = k \, \mathrm{deg} \, Q-(n+1)$.

Let the Jacobian ideal $J(Q)$ be the ideal generated by the partial derivatives $\frac{\partial Q}{\partial
x_i}$.  If we have an element of $\mathcal{H}(X)$ of the form $\frac{K}{Q^{k+1}} \Omega_0$ where $K = \sum_i A_i
\frac{\partial Q}{\partial x_i}$ is a member of the Jacobian ideal, then we may reduce the order of the pole:

\begin{equation}\label{E:reducePoleOrder}
\frac{\Omega_0}{Q^{k+1}} \sum_i A_i \frac{\partial Q}{\partial x_i} = \frac{1}{k} \frac{\Omega_0}{Q^{k}} \sum_i \frac{\partial A_i}{\partial x_i} + \mathrm{exact}\;\mathrm{terms}
\end{equation}

Let $\gamma$ be a cycle in $X$, and let $T(\gamma)$ be a small tubular neighborhood of $\gamma$ in $\mathbb{P}^n
- X$.  Then we may define the {\bf residue map} $\mathrm{Res}: \mathcal{H}(X) \to H^{n-1}(X,\mathbb{C})$ by
\begin{align} \label{residueDef}
    \frac{1}{2\pi i}\int_{T(\gamma)}\frac{P \Omega_0}{Q^k} = \int_{\gamma}\text{Res}(\frac{P \Omega_0}{Q^k})
\end{align}
Let $H$ be the hyperplane class in $H^{n-1}(\mathbb{P}^n, \mathbb{C})$.  We refer to the perpendicular
complement of $\iota^*(H)$ in $H^{n-1}(X,\mathbb{C})$ as the {\bf primitive cohomology} of $X$, and denote it by
$PH(X)$.  The residue map is an isomorphism onto the primitive cohomology. \cite{Griffiths}

Now, consider a family of hypersurfaces $X_{t_1 \cdots t_j}$ given by polynomials $Q_{t_1 \cdots t_j}$, where $t_1, \dots, t_j$ are independent parameters.
We may define a corresponding family
of cycles $\gamma(t_1, \cdots ,t_j)$.  For $(t_1, \dots, t_j)$ in a sufficiently small neighborhood of a fixed
parameter value $(t_1', \dots, t_j')$, $T(\gamma(t_1, \dots, t_j))$ is homologous to $T(\gamma(t_1', \dots, t_j'))$ in $H_n(\mathbb{P}^n-X, \mathbb{C})$.  Thus, we may differentiate as follows:

\begin{align} \label{E:diffUnderIntSign}
    \frac{\partial}{\partial t_i} \int_{T(\gamma(t_1, \dots, t_j))}\frac{P\Omega_0}{Q(t)^k} & = \frac{\partial}{\partial t_i}\int_{T(\gamma(t_1', \dots, t_j'))}\frac{P\Omega_0}{Q(t)^k} \\
    & = -k\int_{T(\gamma(t_1', \dots, t_j'))}\frac{P\Omega_0}{Q(t)^{k+1}}\frac{\partial Q}{\partial t_i} \notag
\end{align}

If $r=\dim_{\mathbb{C}}(H_{n-1}(X))=\dim_{\mathbb{C}}({H^{n-1}(X,\mathbb{C})})$, only $r-1$ derivatives can be
linearly independent. Therefore the periods must satisfy a linear differential equation with coefficients in
$\mathbb{Q}(t_1, \dots, t_j)$ of order at most $r$ --- this is a Picard--Fuchs differential equation.  One may
compute the Picard--Fuchs equation by systematically taking derivatives of $\int_{T(\gamma(t_1, \dots,
t_j))}\frac{P\Omega_0}{Q(t)^k}$ with respect to the various parameters and using \ref{E:reducePoleOrder} to
rewrite the results in terms of a standard basis for $H^{n-1}(X,\mathbb{C})$.  This method is known as the {\bf
Griffiths--Dwork technique}. (See \cite{CK} or \cite{DGJ} for a more detailed discussion.)

\subsubsection{Griffiths--Dwork for the Weierstrass Form} \label{GDforCurves}

Consider the hypersurface
$$Q = y^2 z - 4 x^3 + g_2 x z^2 + g_3 z^3 \ ,$$ the Weierstrass form for a family of elliptic
curves. We illustrate here the Griffiths--Dwork technique, first treating $g_2$ and $g_3$ as independent
parameters. Equation~\ref{E:diffUnderIntSign} tells us that we may differentiate under the integral sign:

\begin{align}
\frac{\partial}{\partial g_2} \int \frac{\Omega_0}{Q} &= -\int \frac{xz^2 \Omega_0}{Q^2} \\
\frac{\partial}{\partial g_3} \int \frac{\Omega_0}{Q} &= -\int \frac{z^3 \Omega_0}{Q^2} \notag
\end{align}

A Groebner basis computation shows that $xz^2$ and $z^3$ are equivalent modulo the Jacobian ideal $J(Q)$.
Using Equation~\ref{E:reducePoleOrder} to reduce the pole order, we find that

\begin{align}
\frac{\partial}{\partial g_2} \int \frac{\Omega_0}{Q} - \frac{\partial}{\partial g_3} \int \frac{\Omega_0}{Q} = \frac{-1}{4 g_2} \int \frac{\Omega_0}{Q}
\end{align}

Now, suppose $g_2$ and $g_3$ are both functions of a single parameter $t$.  We compute:

\begin{align}\label{E:dtWeierstrassIntegral}
\frac{d}{dt} \int \frac{\Omega_0}{Q} & = (\frac{\partial}{\partial g_2} \int \frac{\Omega_0}{Q}) \frac{\partial g_2}{\partial t} + (\frac{\partial}{\partial g_3} \int \frac{\Omega_0}{Q}) \frac{\partial g_3}{\partial t} \\
& = -g_2'(t) \int \frac{xz^2 \Omega_0}{Q^2} - g_3'(t) \int \frac{z^3 \Omega_0}{Q^2} \notag
\end{align}

\begin{align}\label{E:dt2WeierstrassIntegral}
\frac{d^2}{dt^2} \int \frac{\Omega_0}{Q} & = 2 g_2'(t) \int \frac{xz^2 (g_2'(t) xz^2 + g_3'(t) z^3)}{Q^3}\Omega_0 - g_2''(t) \int \frac{xz^2}{Q^2}\Omega \\
& + 2 g_3'(t) \int \frac{z^3 (g_2'(t) xz^2 + g_3'(t) z^3)}{Q^3}\Omega_0 - g_3''(t) \int \frac{z^3}{Q^2}\Omega  \notag \\
& = 2 (g_2'(t))^2 \int \frac{(xz^2)^2}{Q^3} \Omega_0 + 4 g_2'(t) g_3'(t) \int \frac{(xz^2)(z^3)}{Q^3} \Omega_0 \\
& + 2 (g_3'(t))^2 \int \frac{(z^3)^2}{Q^3} \Omega_0 \notag  - g_2''(t) \int \frac{xz^2}{Q^2} \Omega_0 - g_3''(t) \int \frac{z^3}{Q^2} \notag \Omega_0
\end{align}

We may use Equation~\ref{E:reducePoleOrder} together with a Groebner basis computation to rewrite $\frac{d^2}{dt^2} \int \frac{\Omega_0}{Q}$ as a sum of integrals of expressions with $Q^2$ in the denominator:

\begin{align}\label{E:dt2WeierstrassIntegral2}
\frac{d^2}{dt^2} \int \frac{\Omega_0}{Q} & = 2 (g_2'(t))^2 \int \frac{\alpha_1 xz^2 + \beta_1 z^3}{Q^2} \Omega_0 + 4 g_2'(t) g_3'(t) \int \frac{\alpha_2 xz^2 + \beta_2 z^3}{Q^2} \Omega_0 \\
& + 2 (g_3'(t))^2 \int \frac{\alpha_3 xz^2 + \beta_3 z^3}{Q^2} \Omega_0  - g_2''(t) \int \frac{xz^2}{Q^2} \Omega_0 - g_3''(t) \int \frac{z^3}{Q^2} \Omega_0 \notag
\end{align}

Here the $\alpha_j$ and $\beta_j$ are rational functions in $g_2$ and $g_3$.  Note that we have expressed $\frac{d^2}{dt^2} \int \frac{\Omega_0}{Q}$ entirely in terms of $\int \frac{xz^2}{Q^2} \Omega_0 = - \frac{\partial}{\partial g_2} \int \frac{\Omega_0}{Q}$ and $\int \frac{z^3}{Q^2} \Omega_0 = - \frac{\partial}{\partial g_3} \int \frac{\Omega_0}{Q}$.  Since $\frac{d}{dt} \int \frac{\Omega_0}{Q}$ is also written in terms of $\int \frac{xz^2}{Q^2} \Omega_0$ and $\int \frac{z^3}{Q^2} \Omega_0$, we might hope to relate $\frac{d}{dt} \int \frac{\Omega_0}{Q}$ and $\frac{d^2}{dt^2} \int \frac{\Omega_0}{Q}$.  If such a relationship is to exist for an arbitrary choice of $g_2(t)$ and $g_3(t)$, $\int \frac{xz^2}{Q^2} \Omega_0$ and $\int \frac{z^3}{Q^2} \Omega_0$ cannot be independent.  In fact, they are not: $xz^2 \cong  \frac{-3 g_3}{2 g_2}z^3 \; \mathrm{mod} \; J(Q)$, so applying Equation~\ref{E:reducePoleOrder} we find that
\begin{equation}\label{E:dPartialWeierstrassIntegral}
 \int \frac{xz^2}{Q^2} \Omega_0 = \frac{-3 g_3}{2 g_2} \int \frac{z^3}{Q^2} \Omega_0 + \frac{1}{4 g_2} \int \frac{\Omega_0}{Q}.
 \end{equation}

Combining Equations~\ref{E:dtWeierstrassIntegral}, \ref{E:dt2WeierstrassIntegral2}, and \ref{E:dPartialWeierstrassIntegral}, and setting $\Delta = g_2^3 - 27 g_3^2$, we obtain the Picard-Fuchs differential equation for a one-parameter family of elliptic curves in Weierstrass form:

\begin{align}\label{E:tWeierstrassDiffEq}
A_2 \frac{d^2}{dt^2} \int \frac{\Omega_0}{Q} + A_1 \frac{d}{dt} \int \frac{\Omega_0}{Q} + A_0 \int \frac{\Omega_0}{Q} = 0
\end{align}

where

\begin{align}
 A_2 = & 16 \Delta (3 g_2' g_3-2 g_2 g_3')  \\
A_1 = & 16 (9 g_2^2 g_3 (g_2')^2 -
   (7 g_2^3 + 135 g_3^2) g_2' g_3' + 108 g_2 g_3(g_3')^
     2 + \Delta ( -3 g_3 g_2'' +
     2 g_2 g_3''))
  \notag \\
A_0 = &  21 g_2 g_3 (g_2')^3 -
   18 g_2^2 (g_2')^2 g_3' +
   8 g_3'(15 g_2(g_3')^2 -
     \Delta g_2'') -
   4 g_2' (27 g_3(g_3')^2 -
     2 \Delta g_3'') \notag
\end{align}

If we make the substitution $j = g_2^3/\Delta$, then Equation~\ref{E:tWeierstrassDiffEq} reduces to the standard
Picard-Fuchs equation for a one-parameter family of elliptic curves in Weierstrass form, described for example
in \cite{S-H}:

\begin{align}
\frac{d^2}{dt^2} \int \frac{dx}{y} + B_1 \frac{d}{dt} \int \frac{dx}{y} + B_0 \int \frac{dx}{y} = 0
\end{align}

where

\begin{align}
B_1 &= \frac{g_3'}{g_3} - \frac{g_2'}{g_2} + \frac{j'}{j} - \frac{j''}{j'} \\
B_0 &= \frac{(j')^2}{144j(j-1)} + \frac{\Delta'}{12\Delta}\left(B_1 +
\frac{\Delta''}{\Delta'}-\frac{13\Delta'}{12\Delta}\right) \notag
\end{align}

\subsection{Relationship to Toric Geometry} \label{toriccurves}
Another model for the Weierstrass family of elliptic curves comes from toric geometry. This one-parameter family
can be realized as the family of {\bf anticanonical hypersurfaces in the toric variety} $ X = \mathbb{WP}(1,2,3)
$. This family $K$ has defining equation
\[ f_{(\lambda_0, \lambda_1, \lambda_2, \lambda_3)}(x_0,x_1,x_2) = \lambda_0 x_0 x_1 x_2+ \lambda_1 x_1^3+ \lambda_2 x_0^6 + \lambda_3 x_2^2 = 0 \]
in the global homogeneous coordinate ring $ \bb{C}[x_0, x_1, x_2] $ of $X$ (where $ x_0 $, $x_1$, $x_2$ have weights 1, 2, 3 respectively).

The four parameters $ (\lambda_0, \ldots, \lambda_3)$ are redundant. In \cite{CK}, the authors define a
``simplified polynomial moduli space" $ \mathcal{M}_{\mathrm{simp}} $ for $K$ with the property that $
\mathcal{M}_{\mathrm{simp}} $ is a finite-to-one cover (generically) of the actual moduli space $ \mathcal{M} $
for $K$ (which in the present case is just $ \mathbb{WP}(2,3) \simeq \bb{P}^1 $).  Moreover, following
\cite{CK}, we have that $ \mathcal{M}_{\mathrm{simp}} $ is a one-dimensional toric variety, and $ t =
\frac{\lambda_1^2 \lambda_2 \lambda_3}{\lambda_0^4} $ is a coordinate on the torus $ \bb{C}^* \subset
\mathcal{M}_{\mathrm{simp}} $.  We can use the fact that the defining equation is defined only up to a nonzero
constant to set $ \lambda_0 =1 $. Then we can use the natural action of $ T \simeq (\bb{C}^*)^2 $ to set $
\lambda_1 = 4, \lambda_3 = -1 $.  Then our simplified polynomial modulus becomes $ t = -16 \lambda_2 $.

Let $ \phi $ be the rational map $ X \to \bb{P}^2 $ defined by
\[ (x_0, x_1, x_2) \mapsto \left( \frac{x_1}{x_0^2}+\frac{1}{48}, \frac{x_2}{x_0^3}-\frac{x_1}{2 x_0^2}, 1 \right) = \left(x, y, z \right) \]
The image of $ K $ under $ \phi $ has defining equation
\[ y^2 z = 4x^3 - \frac{1}{192} xz^2 + \left(\lambda_2+\frac{1}{13824} \right)z^3 \]
Using the $j$-invariant of an elliptic curve as an affine coordinate on $ \mathcal{M} \simeq \bb{WP}(2,3) $, we have a map $ \mathcal{M}_{\mathrm{simp}} \to \mathcal{M} $ given in affine coordinates by
\[ t \mapsto \frac{1}{1728 t(t-432)} = j \]

The Picard-Fuchs differential equation for $K$ is derived by factoring a Gel'fand-Kapranov-Zelevinsky
hypergeometric equation in \cite[Section 5.1]{LY1}.  The ODE produced there is
\begin{eqnarray} \label{WP123GKZ}
0 & = & \left( \theta^2-12t\left(6 \theta +5 \right) \left( 6 \theta+1 \right) \right) f(t) \\
\nonumber  & = & t \left( t \left( 432t-1 \right) f''(t) + \left(864t-1 \right) f'(t) +60 f(t) \right) \end{eqnarray}
where $ \theta = t \frac{d}{dt} $.  Plugging in $ g_2(t) = 1/192 $, $ g_3(t) = \frac{864t-1}{13824} $ into
Equation~\ref{E:tWeierstrassDiffEq} gives
\begin{equation} 0  =  \frac{1}{393216} \left( t \left( 432t-1 \right) f''(t) + \left(864t-1 \right) f'(t) +60 f(t) \right) \ . \end{equation}
So the methods of GKZ and Griffiths-Dwork yield the same equation in this case (up to an inconsequential overall factor).

\section{K3 Surfaces Polarized by $ {\rm H} \oplus {\rm E}_8 \oplus {\rm E}_8 $} \label{K3PolHE8E8}

\subsection{Normal Form and Moduli Space}
\noindent In \cite{CD1}, the authors presented the full classification of K3 surfaces $({\rm X}, i)$ polarized
by the unimodular even lattice of rank eighteen $ {\rm M} \ = \ {\rm H} \oplus {\rm E}_8 \oplus {\rm E}_8$. The
contents of this section are adapted from \cite{CD1}, where proofs of the results may be found. The crucial
ingredient in this classification is the existence of a {\bf normal form}. That is:
\begin{theorem}[\cite{CD1}]
\label{theo1} Let $({\rm X}, i)$ be an ${\rm M}$-polarized K3 surface. Then, there exists a triple $(a,b,d) \in
\mathbb{C}^3$, with $d \neq 0$ such that $({\rm X}, i)$ is isomorphic to the minimal resolution of the quartic
surface:
\begin{equation}
\label{normalform}
{\rm Q}(a,b,d) \colon \ y^2zw - 4 x^3z + 3 a xzw^2 +b zw^3
- \frac{1}{2} ( d z^2 w^2 + w^4 ) \ = \ 0. \
\end{equation}
\end{theorem}
\noindent The equation in Theorem $\ref{theo1}$ can be regarded as the analog, for ${\rm M}$-polarized K3
surfaces, of the Weierstrass form in the classical theory of elliptic curves. It is a compactification of
Inose's $(a,b)$ family from Equation \ref{InoseIntro}.
\par Two distinct quartics as in $(\ref{normalform})$ may correspond to isomorphic polarized K3 surfaces.
\begin{theorem}[\cite{CD1}]
\label{theo2} Two quartics $ {\rm Q}(a_1, b_1, d_1) $ and ${\rm Q}(a_2, b_2, d_2)$ determine isomorphic ${\rm
M}$-polarized K3 surfaces as their minimal resolutions if and only if:
\[ (a_2, \ b_2, \ d_2) \ = \ (\lambda^2 a_1, \ \lambda^3 b_1, \ \lambda^6 d_1) \]
for some parameter $\lambda \in \mathbb{C}^*$.
\end{theorem}
\noindent One obtains therefore a {\bf coarse moduli space} for ${\rm M}$-polarized K3 surfaces in the form of the open variety:
\begin{equation}
\label{moduliabd}
\mathcal{M}_{{\rm M}} \ = \ \left \{ \ [a,b,d] \ \in \mathbb{WP}(2,3,6) \ \middle \vert \ d \neq 0 \ \right \}
\end{equation}
as well as a pair of {\bf fundamental $\mathcal{W}$-invariants} $({\mathcal W}_1, {\mathcal W}_2)$:
\[ \mathcal{W}_1 \ = \ \frac{a^3}{d}, \ \  \mathcal{W}_2 \ = \ \frac{b^2}{d}  .\]
\par We turn to Hodge theory. Denote by ${\rm T}$ the rank-four even indefinite lattice ${\rm H} \oplus {\rm H}$. Then, following the general
framework of \cite{dolga1}, one notes that the
Hodge structure of an ${\rm M}$-polarized K3 surface $({\rm X},i)$ can be seen as a point in the classifying space:
\begin{equation}
\label{hodgeperiodspace}
\Gamma \backslash \Omega
\end{equation}
where
$$ \Omega \ = \ \left \{ \
\omega \in \mathbb{P}({\rm T} \otimes \mathbb{C} ) \ \middle \vert \ (\omega, \omega)=0, \ (\omega, \overline{\omega}) >0 \
 \right \} $$
 and $\Gamma$ is the group of integral isometries of the lattice ${\rm T}$, with its standard action on $\Omega$. Upon further inspection,
 the quotient $(\ref{hodgeperiodspace})$ turns out to be identical to the {\bf classical Hilbert modular surface}:
\begin{equation}
\mathcal{H}_1 \ = \
\left ( {\rm PSL}(2, \mathbb{Z}) \times {\rm PSL}(2, \mathbb{Z}) \right ) \rtimes \mathbb{Z}/2 \mathbb{Z}
\  \backslash \ \mathbb{H} \times \mathbb{H}
\end{equation}
where $\mathbb{H}$ denotes a copy of the complex upper-half plane. This reinterpretation reflects, at a Hodge theoretic level,
the geometric one-to-one correspondence existing between ${\rm M}$-polarized K3 surfaces and abelian surfaces that split as a product
of two elliptic curves.
\par The Global Torelli Theorem (see \cite{dolga1}) asserts then that the period map:
$$ {\rm per} \colon \mathcal{M}_{{\rm M}} \rightarrow \mathcal{H}_1 $$
is an analytic isomorphism. This isomorphism can be made completely explicit. Let $j \colon \mathbb{H} \rightarrow \mathbb{C}$ be
the classical elliptic modular function. There are two important modular functions on the Hilbert surface
$\mathcal{H}_1$, given, on the cover $\mathbb{H} \times \mathbb{H}$,  by the two symmetric functions:
\[ \sigma(\tau_1, \tau_2) \ = \ j(\tau_1) + j(\tau_2), \ \  \pi(\tau_1, \tau_2) \ = \ j(\tau_1) \cdot j(\tau_2). \]
They can be seen naturally as coordinates on the surface $\mathcal{H}_1$ realizing an identification:
\[ (\sigma, \pi) \colon \mathcal{H}_1 \ \stackrel{\simeq}{\longrightarrow} \ \mathbb{C}^2. \]
\begin{theorem}[\cite{CD1}]
\label{theo3} The inverse period map ${\mathrm per}^{-1} \colon \mathcal{H}_1 \rightarrow \mathcal{M}_{{\rm M}}
$ is given by:
$$ {\mathrm per}^{-1} \ = \ \left [ \pi ^{1/3}, \ ( \pi - \sigma +1)^{1/2}, \ 1 \right ]. $$
In other words, the $\mathcal{W}$-invariants of an ${\rm M}$-polarized K3 surface $({\rm X},i)$ are linked to the
periods of ${\rm X}$ by the formulas:
$$ \mathcal{W}_1 = \pi, \ \ \mathcal{W}_2 = \pi - \sigma + 1. $$
\end{theorem}

\subsection{Loci Supporting Lattice Enhancements} \label{enhancement}
Since the rank 18 lattice ${\rm M}$ is unimodular,  the only possible enhancements to ${\rm M}'$-polarized
lattices with ${\rm M} \subset {\rm M}'$ and ${\rm M}'$ of rank 19 are those for which ${\rm M}' = {\rm M}_n =
{\rm M} \oplus \langle - 2 n \rangle = {\rm H} \oplus {\rm E}_8 \oplus {\rm E}_8 \oplus \langle - 2 n \rangle$
for some $n \in {\mathbb{N}}$. The $ {\rm M}_n $ polarized K3 surfaces for $n>1$ correspond, via the
Shioda-Inose construction, to products of elliptic curves with an $n$-isogeny between them. The extra algebraic
cycle on the K3 corresponds to the graph of the $n$-isogeny on the product of the two elliptic curves.
Determining the subloci of $ \bb{WP}(2,3,6) $ on which these enhancements occur thus reduces to the problem of
finding relations between the $j$-invariants of pairs of elliptic curves. This, of course, is a classical
problem with a rich history.

Consider the map $ \phi_n : \mathbb{H} \to \mathbb{P}^2 $ given by $ \tau \mapsto [j(\tau), j(n \tau), 1] = [x,
y, z] $, and denote by $ X_0(n) $ the closure of the image of $ \phi_n $.   $ \phi_n(\tau) = \phi_n(\tau') $ if
and only if
\[ \tau' = \frac{a \tau + b}{c \tau + d} \ , \ n \tau' = \frac{\alpha n \tau + \beta}{\gamma n \tau +
\delta} \]
for some
\[ \left( \begin{array}{cc} a & b \\ c & d \end{array} \right), \: \left( \begin{array}{cc}
\alpha & \beta \\ \gamma & \delta \end{array} \right) \in SL_2(\bb{Z}) \]
Examining these equations, we see that
\[\left( \begin{array}{cc} na & nb \\ c & d \end{array} \right)=\pm \left( \begin{array}{cc} n \alpha & \beta \\
n \gamma & \delta \end{array} \right) \] which is possible if and only if $ \tau' = \frac{a \tau + b}{c \tau +d}
$, for \[ \left( \begin{array}{cc} a & b \\ c & d \end{array} \right) \in \Gamma_0(n) \ . \] Thus $ X_0(n) $ is
a compactification of $ \Gamma_0(n) \backslash \mathbb{H} $, and parametrizes (ordered) pairs of $n$-isogenous
elliptic curves.

The genus of $ X_0(n) $ for many small $n$ can be found in \cite{Cohn}.  The equation $ \Phi_n(x,y) = 0 $ for $
X_0(n) $ on the affine patch $ z= 1 $ is the classical modular equation for $ \Gamma_0(n) $\footnote{The
$j$-invariant we use is normalized so that $ j(i)=1 $, not 1728 as in many works on modular equations in the
literature.}.  Because it is symmetric in $x$ and $y$, it can be written as a polynomial in terms of the
elementary symmetric functions $ \pi = xy$, $ \sigma = x + y $.  Let $ \Phi_n^{+n}(\pi, \sigma) $ be the
corresponding polynomial such that $ \Phi_n^{+n}(xy, x+y) = \Phi_n(x,y) $.  Note that
$$ (j(\tau) j(n \tau), j(\tau) + j(n \tau)) = (j(\tau') j(n \tau'), j(\tau') + j(n \tau')) $$
if and only if $ (j(\tau), j(n \tau)) = (j(\tau'), j(n \tau') $ or $ (j(n \tau'), j(\tau')) $,
which will occur if and only if $ \tau' = \frac{a \tau +b}{c \tau + d} $ for some
\[ \left( \begin{array}{cc} a & b \\ c & d \end{array} \right) \in \Gamma_0(n)+n \ . \]
We may thus view $ \Phi_n^{+n}(\pi,\sigma) = 0 $ as the modular relation for $ \Gamma_0(n)+n $. Then, using the
relationships $ \pi = a^3, \; \sigma = a^3 - b^2 + 1 $ from \cite{CD1}, we can rewrite $ \Phi_n(x,y) $ as a
polynomial $ \Psi_n^{+n} (a^3, b^2) $, and view $ \Psi_n^{+n} (a^3, b^2) = 0 $ as the equation for a curve $
Y_0(n)+n \subset \bb{WP}(2,3,6) $ on the affine patch $ d = 1 $ of the weighted projective plane with
homogeneous coordinates $ [a , b , d ] $ of weights 2, 3, and 6 respectively.  Note $ X_0(n)+n \simeq Y_0(n)+n
$.  $ Y_0(n)+n $ is exactly the moduli space of $ {\rm M}_n $-polarized K3 surfaces as a submoduli space of the
moduli space of $ {\rm M}$-polarized K3 surfaces.

Thus the theory of modular curves and modular equations can be used to analyze the moduli of ${\rm
M}_n$-polarized K3 surfaces.  For example, the genus of $ Y_0(n)+n $ can be found in the literature, and in the
case where $ Y_0(n)+n $ has genus zero, it can be parametrized by the hauptmodul for $ \Gamma_0(n)+n $.

\subsubsection{Examples: $n = 2 $, $3$, and $6$}
Since $ \Gamma_0(6) $ is contained in $ \Gamma_0(3) $ and $ \Gamma_0(2) $, it is natural to work over $ X_0(6) $. Choosing as affine coordinate on $ X_0(6) \simeq \bb{P}^1 $ the hauptmodul $ t $ for $ \Gamma_0(6) $, the family of elliptic curves over $ X_0(6) $ was shown in \cite{Beauville} to form a rational elliptic surface $ S $ with four singular fibers, of Kodaira types $ I_1, I_2, I_3, I_6 $.  Beauville gave a model for this surface: \[ S = \{ ([x, y, z], t) \in \bb{P}^2 \times \bb{P}^1 : (x+y)(y+z)(z+x) + t xyz = 0 \} \]
with the singular fibers of types $ I_1, I_2, I_3, I_6 $ over the points \[ t=-8,1,0,\infty \] respectively.  Another model for this surface is given in \cite{Hadano}: \[ S' = \{ ([x, y, z],[\alpha, \beta]) \in \bb{P}^2 \times \bb{P}^1 : y^2z-2(\alpha+\beta)xyz+2\alpha \beta^2yz^2+x^3=0 \} \]  An explicit isomorphism between $ S $ and $ S' $ is \[
\xymatrix{
S \ar[d] \ar[r]^f & S' \ar[d] \\
\bb{P}^1 \ar[r]^g & \bb{P}^1 }
\] where $ f $ is defined by\[ ([x, \ y, \ z],t) \mapsto ([tyz , t(z^2+(2+t)yz+y^2+x(z+y)) , -z(z+x)], [t/2 , 1]) \] and $ g $ is defined by \[ t \mapsto [t/2 , 1] \]

To describe the $2$, $3$, and $6$-isogenies explicitly, it is more convenient to use $ S' $.  We will denote by
$ S'_{[\alpha, \ \beta]} $ the fiber of $ S'$ over the point $ [\alpha, \ \beta] \in \bb{P}^1 $.  From
\cite{Hadano}, we have that $ S'_{[\alpha, \ \beta]} $ is $2$-isogenous to $ S'_{[-2\beta,\ \alpha]} $,
$3$-isogenous to $ S'_{[-\alpha - 4\beta, \ -2\alpha +\beta]}$, and $6$-isogenous to $ S'_{[4\alpha - 2\beta, \
-\alpha - 4\beta]} $.  The isogenies are not given in \cite{Hadano}, but \cite{Cassels} gives an explicit
2-isogeny, and \cite{Kuwata} gives an explicit 3-isogeny for families of elliptic curves.  Comparing these
families with $ S' $ and noting that the 6-isogeny $ S'_{[\alpha, \ \beta]} \to S'_{[4\alpha - 2\beta, -\alpha -
4\beta]} $ is the composition of the 2-isogeny $ S'_{[\alpha:\beta]} \to S'_{[-2\beta,\alpha]} $ with the
3-isogeny $ S'_{[-2\beta,\alpha]} \to S'_{[4\alpha - 2\beta, -\alpha - 4\beta]} $, we can derive explicitly what
the isogenies $ \phi_n: S' \to S' $ are for $ n = 2,3, 6 $.  Denote $ \phi_n(([x,y,z],[\alpha,\beta]) $ by $
([x'_n, y'_n,z'_n], [\alpha'_n, \beta'_n]) $. Then $ x'_n, y'_n, z'_n,\alpha'_n, \beta'_n $ are given below for
$n=2,3$.  (In the interest of space, we omit the corresponding data for $n=6$, but we note that it can easily be
deduced, since $ \phi_6 $ is the composition of $ \phi_2 $ and $ \phi_3 $.)
\begin{eqnarray*}
x'_2 & = & x(\beta^2z-x)(2\alpha \beta z - x) \\
y'_2 & = & (\beta-2\alpha)x^3 + 2\alpha \beta^3yz^2 - 2\beta^2xyz + x^2y \\
z'_2 & = & (\beta^2z-x)^2z \\
\alpha'_2 & = & -2\beta \\
\beta'_2 & = & \alpha \\
& & \\
x'_3 & = & x ((4\beta^2z - 3x)x^2 + 4\alpha \beta xz(-3\beta^2z + 2x) + 4\alpha^2z(3\beta^4z^2 - 3\beta^2xz + x^2)) \\
y'_3 & = & -8\alpha^3z(3\beta^6z^3 - 3\beta^2x^2z + x^3) +
 x^3(-8\beta^3z + 12\beta x - 3y) - \\
 & & 12\alpha^2\beta z(2\beta^2z - x)
  (-2x^2 + \beta yz) + 6\alpha x(4\beta^4xz^2 - 9\beta^2x^2z + 2x^3 +
   2\beta ^3yz^2) \\
   & & - 6 \zeta_3 (-4\alpha \beta^3xz^2 + x^3 + \alpha^2z(8\beta^4z - 4\beta^2x))
 (\alpha(\beta^2z - x) - \beta x + y) \\
 z'_3 & = & x^3z \\
 \alpha'_3 & = & -\alpha-4\beta  \\
 \beta'_3 & = & -2\alpha+\beta \end{eqnarray*}
where $ \zeta_3 $ is a primitive third root of unity.

In order to understand $ Y_0(n)+n $ for $n=2,3,6$ as a sublocus of $ \bb{WP}[2,3,6] $, there are two ways one
could proceed.  One is to find the defining equation $ \Psi_n^{+n}(a,b,d) = 0 $ for $ Y_0(n)+n $, and another
--- since in these cases $ Y_0(n)+n $ is genus zero --- is to find a parametrization $[a(t), b(t), d(t)]$ for the
curve.  Of course, these two points of view are closely related.  Given a parametrization one can eliminate the
parameter to obtain a defining equation, and conversely given the defining equation for a genus zero curve one
can employ standard algorithms (for example in \cite{AB}) to obtain a parametrization.

Defining equations for $ X_0(n) $ for small $n$ are readily available (for example, in MAGMA), and from these we can obtain defining equations $ \Psi_n^{+n}(a,b,d) = 0 $ for $ Y_0(n)+n $.  We give $ \Psi_n^{+n}(a,b,d) $ for $ n=2,3 $ here. $ \Psi_6^{+6}(a,b,d) $ is a degree 72 polynomial with large coefficients, so we will omit it.
\begin{eqnarray*}
\Psi_2^{+2}(a,b,d) & = & 64 a^9-192 b^2 a^6-21360 d a^6+192 b^4 a^3+1792857 d^2
   a^3 \\
   & & -83424 b^2 d a^3 -64 b^6-1771561 d^3-175692 b^2
   d^2-5808 b^4 d \\
\Psi_3^{+3}(a,b,d) & = & 729 a^{12}-2916 b^2 a^9-129551076 d a^9+4374 b^4
   a^6 \\
   & & +5754777529878 d^2 a^6 -1298340252 b^2 d a^6-2916
   b^6 a^3 \\
   & & +77703185570076 d^3 a^3- 3974452231068 b^2 d^2
   a^3 \\
   & & -733336092 b^4 d a^3+729 b^8 +262254607552729
   d^4 \\
   & & -262365230658916 b^2 d^3 +110638660374 b^4 d^2-15554916 b^6 d
\end{eqnarray*}

Given a parametrization $ (x(t_n), y(t_n)) $ of $ X_0(n) $ on the affine patch $ z \neq 0 $, using coordinates $(b^2,d)$ on the affine patch $ a \neq 0 $ of $ \bb{WP}(2,3,6)$, we can construct a parametrization of $ Y_0(n)+n $ by setting \begin{eqnarray*}
b^2(t_n) & = &  \frac{(x(t_n)-1)(y(t_n)-1)}{x(t_n)y(t_n)} \\
d(t_n) & = & \frac{1}{x(t_n) y(t_n)} \end{eqnarray*}
\cite{Maier} gives such a parametrization $ (x(t_n), y(t_n)) $ by the hauptmodul $ t_n $ for $ \Gamma_0(n) $ for $ n \leq 25 $ such that $ X_0(n) $ is a genus 0 curve.  These parametrizations for $ n=2,3,6 $ give \begin{eqnarray*}
b^2(t_2) & = & \frac{(t_2-512)^2 (t_2-8)^2 (t_2+64)^2}{(t_2+16)^3 (t_2+256)^3} \\
d(t_2) & = & \frac{2^{12} 3^6 t^3}{(t+16)^3 (t_2+256)^3} \\
& & \\
b^2(t_3) & = & \frac{\left(t_3^2-486 t_3-19683\right)^2 \left(t_3^2+18
   t_3-27\right)^2}{(t_3+27)^2 \left(t_3^2+246 t_3+729\right)^3} \\
d(t_3) & = & \frac{2^{12} 3^6 t_3^4}{(t_3+3)^3 (t_3+27)^2 (t_3+243)^3} \\
& & \\
b^2(t_6) & = & \left(t_6^2+12 t_6+24\right)^2 \left(t_6^2+36
   t_6+216\right)^2 \times \\
   & & \left(-t_6^4+504 t_6^3+13824 t_6^2+124416
   t_6+373248\right)^2 \times \\
   & & \frac{ \left(t_6^4+24 t_6^3+192 t_6^2+504
   t_6-72\right)^2}{\left( (t_6+6) (t_6+12)
   (t_6^3+18 t_6^2+84 t_6+24) (t_6^3+252
   t_6^2+2^4 3^5 t_6+2^6 3^5) \right)^3} \\
d(t_6) & = & \frac{2985984 t_6^7 (t_6+8)^5 (t_6+9)^5}{\left( (t_6+6) (t_6+12)
   (t_6^3+18 t_6^2+84 t_6+24) (t_6^3+252
   t_6^2+2^4 3^5 t_6+2^6 3^5) \right)^3}
\end{eqnarray*}

Another interesting parametrization of $ Y_0(n)+n $ can be derived from the ``two-valued" modular equations derived, following Fricke, in \cite{Cohn}.

\subsection{Picard-Fuchs Differential Equations}

Let us set $ a = 1 $ in Equation \ref{normalform} and consider the resulting polynomial $Q = y^2zw - 4 x^3z + 3xzw^2 + bzw^3 - \frac{1}{2} ( d z^2 w^2 + w^4)$.  (We have simply reduced to the affine patch $ a \neq 0 $ of the parameter space $\bb{WP}(2,3,6)$.)  Applying the Griffiths--Dwork technique to $\int \frac{\Omega_0}{Q}$ yields a pair of second-order Picard--Fuchs equations:

\begin{align}  \label{GDbd1}
\frac{\partial^2}{\partial b^2} \int \frac{\Omega_0}{Q} - 4 (d \frac{\partial^2}{\partial d^2} \int \frac{\Omega_0}{Q} + \frac{\partial}{\partial d} \int \frac{\Omega_0}{Q}) &= 0\\
\label{GDbd2}
(-1 + b^2 + d)\frac{\partial^2}{\partial b^2} \int \frac{\Omega_0}{Q} + 2b\frac{\partial}{\partial b} \int \frac{\Omega_0}{Q} + 4bd\frac{\partial^2}{\partial bd} \int \frac{\Omega_0}{Q} & \\
+ 2d\frac{\partial}{\partial d} \int \frac{\Omega_0}{Q} + \frac{5}{36} \int \frac{\Omega_0}{Q} &=0 \notag
\end{align}

We can use the relationship between $b, d $ and the $j$-invariants of elliptic curves from Theorem \ref{theo3}
to write $ b^2 = \frac{(j_1-1)(j_2-1)}{j_1 j_2 } $ and $ d = \frac{1}{j_1 j_2} $. Here $j_1$ and $j_2$ are the
$j$-invariants of the two elliptic curves ${\rm E}_1$ and ${\rm E_2}$ whose product corresponds to $X(1, b, d)$.
  Let $ {\rm E}_i $ have affine Weierstrass model
\[ y^2 = 4 x^3 - g_2^{(i)} x - g_3^{(i)} \] for $ i = 1, 2$. Then we can rewrite Equations \ref{GDbd1} and \ref{GDbd2} in terms of $j_1$ and $j_2 $.  The resulting system decouples (no mixed partials appear).
The system reduces (after taking appropriate linear combinations of the resulting equations) to
\begin{eqnarray*}
 \ \ \ 0 & = & 72 j_1 \left((2 j_1-1)
   F^{(1,0)}(j_1,j_2)+2 (j_1-1)
   j_1 F^{(2,0)}(j_1,j_2)\right)-5
   F(j_1,j_2) \\
 \ \ \ 0 & = & 72 j_2 \left((2
   j_2-1) F^{(0,1)}(j_1,j_2)+2
   (j_2-1) j_2
   F^{(0,2)}(j_1,j_2)\right)-5
   F(j_1,j_2)
   \end{eqnarray*}
where $ F^{(i,j)}(j_1, j_2) = \frac{\partial^{i+j} F}{\partial j_1^i \partial j_2^j} $.

To solve this system, one need merely solve each ODE separately, then take products of the solutions.  Each of
these ODEs separately is a Picard-Fuchs differential equation satisfied by the periods of the form $
\omega^{(i)} = \left( g_2^{(i)} \right)^{1/4} \frac{dx}{y} $. Thus periods satisfying the Picard-Fuchs system
arising via Griffiths-Dwork are simply products of periods of $ \omega^{(1)} $ and $ \omega^{(2)} $ (c.f.
\cite[Theorem 1.1]{LY2}).

Now, consider a one-parameter family $ \mathcal{F} $ of ${\rm M}$-polarized K3 surfaces obtained by treating $b$
and $d$ as functions of a single parameter $t$.  We may use the Griffiths--Dwork technique to analyze this
family, just as we computed the Picard--Fuchs equation for a one-parameter family of elliptic curves in
Section~\ref{GDforCurves}.  The result is generically a fourth-order ODE, which we do not reproduce in full
here. The Picard--Fuchs equation for $ \mathcal{F} $ will reduce to a third-order ODE precisely when $
\mathcal{F} $ is an ${\rm M}_n$-polarized family.

Let $ j_1(t) $, $ j_2(t) $ be two functions of a complex variable $t$ such that $ j_1(t) + j_2(t) $ and $ j_1(t)
j_2(t) $ are rational functions of $t$.  In this case $ b^2(t) = \frac{(j_1(t)-1)(j_2(t)-1)}{j_1(t) j_2(t) } $
and $ d(t) = \frac{1}{j_1(t) j_2(t)} $ are also rational functions of $t$, and we may write the Picard--Fuchs
equation for $ \mathcal{F} $ in terms of $ j_1(t) $ and $ j_2(t) $. The coefficient $ r_4(t) $ of $
\frac{\mathrm{d}^4}{\mathrm{d} t^4} \int \frac{\Omega_0}{Q}$ in the Picard--Fuchs ODE then becomes
\begin{equation*}
 144 ((j_1(t)-1)(j_2(t)-1))^3  (j_1(t) j_2(t))^4 (j_1(t)-j_2(t))^7  (j_1'(t) j_2'(t))^2 \left( \Box(j_2(t))-\Box(j_1(t)) \right) 
\end{equation*}
where $$ \Box(j(t)) = j'(t)^2 \frac{36 j(t)^2-41 j(t)+32}{144(j(t)-1)^2 j(t)^2}+ \frac{1}{2} \{ j(t), t\} $$ and
$$ \{j(t), t\} = \frac{2 j'(t) j'''(t)-3 j''(t)^2}{2 j'(t)^2}$$ is the Schwarzian derivative.

If $ j_1(t) $ and $ j_2(t) $ are both nonconstant, then $ r_4(t) $ will vanish if and only if either $ j_1(t) =
j_2(t) $ (in which case the family of K3 surfaces is ${\rm M}_1$-polarized) or $ \Box(j_1(t)) = \Box(j_2(t)) $.
This observation motivates the following theorem.

\begin{theorem}
\label{boxtheorem} Let $ j_1(t), j_2(t) $ be nonconstant functions of a complex variable $t$ such that $ j_1(t)
+ j_2(t) $ and $ j_1(t) j_1(t) $ are rational functions of $t$. Then \begin{equation} \label{mastereqn}
\Box(j_1(t)) = \Box(j_2(t)) \end{equation} if and only if $ (j_1(t), j_2(t)) $ is a parametrization of $ X_0(n)
$ for some $n \geq 1$.
\end{theorem}

\begin{proof}
The Picard-Fuchs ODE for $ \mathcal{F} $, suitably normalized, is the tensor product of the Picard-Fuchs ODE's
of the two pencils of elliptic curves over $ \bb{P}^1_t $ with functional invariants $ j_1(t), j_2(t) $
respectively. If these second-order ODE's $ L_1 =0 , L_2 = 0 $ are in projective normal form
\begin{eqnarray*}
L_1 & = & \dd{f}{t}{2} + p_2(t) f  \\
L_2 & = & \dd{g}{t}{2} + q_2(t) g
\end{eqnarray*}
then $ p_2(t) = \Box(j_1(t)) $ and $ q_2(t) = \Box(j_2(t)) $.  Their tensor product is

\begin{eqnarray} \nonumber 0 & = & H^{(4)}(t)+\frac{q_2'(t)-p_2'(t)}{p_2(t)-q_2(t)} H'''(t) + 2 (p_2(t)+q_2(t)) H''(t)+ \\
& & \frac{p_2(t) \left(p_2'(t)+5 q_2'(t)\right)-q_2(t) \left(5
   p_2'(t)+q_2'(t)\right)}{p_2(t)-q_2(t)}H'(t)+ \\
\nonumber   & &
  \left( (p_2(t)-q_2(t))^2+p_2''(t)+q_2''(t)+\frac{q_2'(t)^2-p_2'(t)^2}{p_2(t)-q_2(t)} \right) H(t)  \end{eqnarray}

According to \cite{Fano}, this fourth-order equation factorizes as a third-order equation times a first-order
equation if and only if $ p_2(t) = q_2(t) $, i.e. if and only if $ \Box(j_1(t)) = \Box(j_2(t)) $.  On the other
hand,  the Picard-Fuchs equation of $ \mathcal{F} $ has third order if and only if $ \mathcal{F} $ is $ {\rm
M}_n $-polarized, and this occurs if and only if the two pencils of elliptic curves are fiberwise
$n$-isogenous--or in other words, if and only if $ (j_1(t), j_2(t)) \in X_0(n) $ for all $t$.
\end{proof}
For this reason, we call Equation \ref{mastereqn} {\bf the master equation for modular parametrization} of
modular equations for the elliptic modular function $j(\tau)$.

\begin{corollary} If two pencils of elliptic curves over $ \bb{P}^1_t $ admit a fiberwise $n$-isogeny for $n$ such that $ X_0(n)+n $ is genus 0, then the projective normal forms of their Picard-Fuchs differential operators are identical. \end{corollary}

\subsubsection{Modular Relations and Differential Identities for Hauptmoduls} \label{hauptmoduls}
Theorem \ref{boxtheorem}, which relates the existence of a parametrized modular relation between $ j(\tau) $ and
$ j( n \tau) $ to a differential identity involving a parametrization, can be generalized to other hauptmoduls.
If $ h $ is a hauptmodul for a genus 0 modular group $ \Gamma $, then $h$ will satisfy a Schwarzian differential
equation

\begin{eqnarray} \label{SchwarzianDE} h'(\tau)^2 Q_\Gamma(h(\tau)) + \frac{1}{2} \{ h(\tau), \tau) \} & = & 0 \end{eqnarray}
where $ Q_\Gamma(h) $ is a rational function we shall call the ``Q-value for h."  A list of ``Q-values" for (suitably normalized) hauptmoduls for all genus 0 modular groups are given in \cite{LW}.

As was noted in \cite{HM,Har}, if the hauptmodul $h$ for $ \Gamma $ can be expressed as a function $ h_1 (t) $
of a hauptmodul $t$ for a genus 0 group $ \Gamma' $, then
\begin{eqnarray} \label{onebox} h_1'(t)^2 Q_\Gamma(h_1(t)) + \frac{1}{2} \{ h_1(t), t \} & = & Q_{\Gamma'}(t) \end{eqnarray}
This identity follows quite easily from formally writing $ h_1(t) = h(\tau(t)) $, and applying the chain rule to
Equation \ref{SchwarzianDE}.  If we can also write $ h( n \tau) $ as a function $ h_2(t) $, then we have an
analogous result to Equation \ref{onebox} with $ h_1(t) $ replaced by $ h_2(t) $, and hence we have

\begin{theorem} \label{boxgeneral} If $h$ and $t$ are hauptmoduls for genus zero modular groups $ \Gamma $ and $ \Gamma'$ respectively, and if $ h(\tau), h( n \tau) $ can both be expressed as rational functions of $t$ then
\begin{eqnarray}
h_1'(t)^2 Q_\Gamma(h_1(t)) + \frac{1}{2} \{ h_1(t), t \} & = & h_2'(t)^2 Q_\Gamma(h_2(t)) + \frac{1}{2} \{ h_2(t), t \} \\
\nonumber & = & Q_{\Gamma'}(t) \end{eqnarray}
\end{theorem}

Such a situation will occur when $ h(\tau), h(n \tau) $ satisfy a modular equation of genus zero.  Theorem \ref{boxtheorem} covers the special case $ \Gamma = PSL(2, \bb{Z}) $, $\Gamma' = \Gamma_0(n) $.

For example, \cite{CM} discusses the hauptmoduls satisfying modular equations of levels two and three, and gives formulae for the modular equations.  Parametrizing the curves defined by these equations, computing the left-hand side of Equation \ref{onebox}, and comparing with list of ``Q-values" in \cite{LW}, we can verify Theorem \ref{boxgeneral} in these cases (and identify the group $ \Gamma' $ if it is not already known).  We illustrate this for $ \Gamma = \Gamma_0(3)+3 $ below.

The hauptmodul $ h $ for $ \Gamma = \Gamma_0(3)+3 $ has Q-value $ \frac{h^2-48 h+7560}{4 \left(h^2-24 h-2772\right)^2} $ and satisfies the level two modular equation $ \Phi_2(h(\tau), h( 2 \tau)) = 0 $ where
\[ \Phi_2(h_1,h_2) = h_1^3-h_2^2 h_1^2+17343 h_2 h_1+h_2^3+741474 (h_1+h_2)+1566 (h_1^2+h_2^2)+28166076 \]
The curve defined by $ \Phi_2(h_1,h_2) = 0 $ can be parametrized by setting
\[ h_1(t) = \frac{-512 t^3+804 t^2-12 t+1}{(1-6 t)^2 t}, \> h_2(t) = \frac{244 t^3-30 t^2-6 t+1}{t^2-6 t^3} \]
Evaluating the left-hand side of \ref{onebox}, we see that \[ Q_{\Gamma'}(t) = \frac{2848 t^4-800 t^3+108 t^2+4 t+1}{4 t^2 \left(120 t^3-68 t^2+2 t+1\right)^2} \] and by comparing with \cite{LW} we see that $ t $ is a hauptmodul for $ \Gamma' = \Gamma_0(6)+3 $.

\subsection{Relationship to Toric Geometry}
Another model for the two parameter family of ${\rm M}$-polarized K3 surfaces comes from toric geometry.  The
2-parameter family can be realized as the family of anticanonical hypersurfaces in the mirror (polar) toric
variety $ X $ to $ \bb{WP}(1,1,4,6) $. The family of anticanonical K3's $K$ has defining equation
$$
f_{(\lambda_0, \lambda_1, \ldots, \lambda_5)} (x_0,x_1,x_2,x_3) = \lambda_0 x_0 x_1 x_2 x_3 + \lambda_1 x_0^{12}
+ \lambda_2 x_1^{12} + \lambda_3 x_2^3+\lambda_4 x_3^2 + \lambda_5 x_0^6 x_1^6 = 0 $$ in the global homogeneous
coordinate ring $ \bb{C}[x_0,x_1,x_2,x_3] $ of $X$ (where $ x_0,x_1,x_2,x_3 $ have weights $1,1,4,6$
respectively).

As in Section \ref{toriccurves}, these six parameters $ (\lambda_0, \lambda_1, \ldots, \lambda_5) $ are
redundant.  Let $ \mathcal{M}_{\mathrm{simp}} $ be the simplified polynomial moduli space for $ K $. $
\mathcal{M}_{\mathrm{simp}} $ is a two-dimensional toric variety and is a finite-to-one cover (generically) of
the actual moduli space $ \mathcal{M} $ (which in the present case is $ \bb{WP}(2,3,6) $).  We will use $ (z_1,
z_2) = (\frac{\lambda_3^2 \lambda_4^3 \lambda_5}{\lambda_0^6}, \frac{\lambda_1 \lambda_2}{\lambda_5^2}) $ as
affine coordinates on the torus $ \left( \bb{C}^* \right)^2 \subset \mathcal{M}_{\mathrm{simp}} $.  We can use
the fact that the defining equation is only defined up to an overall nonzero constant to set $ \lambda_4 = 1 $.
Then we can use the natural action of $ T \simeq (\bb{C}^*)^3 $ on $X$ to set $ \lambda_1=\lambda_2=-1/2 $, and
$\lambda_3 = -4 $.  Then our simplified polynomial moduli are $ (z_1, z_2) = (\frac{16 \lambda_5}{\lambda_0^6},
\frac{64}{\lambda_5}) $.

Let $ \phi $ be the rational map $ X \to \bb{P}^3 $ defined by
\[ (x_0, x_1, x_2, x_3) \mapsto
\left(\frac{x_2}{x_0^2 x_1^2} - \frac{\lambda_0^2}{48}, \frac{x_3}{x_0^3 x_1^3}+\frac{\lambda_0 x_2}{2 x_0^2
x_1^2}, \frac{x_1^6}{x_0^6}, 1\right) = (x, y, z, w) \]

The image of $ K $ under $ \phi $ has defining equation \begin{eqnarray*}
0 & = & y^2 z w - 4x^3z+ \frac{\lambda_0^4}{192} xzw -\frac{1}{2} \left( z^2w^2+w^4 \right) + \left( \frac{-\lambda_0^6}{13824} + \lambda_5 \right) z w^3 \\
& = & y^2 z w - 4x^3z+ \frac{\lambda_0^4}{192} xzw -\frac{1}{2} \left( z^2w^2+w^4 \right) + \left( \frac{-\lambda_0^6}{13824} + \lambda_5 \right) z w^3
\end{eqnarray*}
Comparing this with the Inose normal form for ${\rm M}$-polarized K3 surfaces, we get a map $
\mathcal{M}_{\mathrm{simp}} \to \bb{WP}(2,3,6) $ given in affine coordinates by
\[ (z_1, z_2) \mapsto (\frac{1}{144^3 z_1^2 z_2}, \frac{864 z_1-1}{144^3 z_1^2 z_2}) = (a^3, b^2) \]
The Picard-Fuchs equations for this family has been computed from GKZ methods (e.g. by Lian-Yau in \cite{LY2})
as
\begin{eqnarray}
0 & = & L_1 f(z_1, z_2) \\
\nonumber  & = & \left( \theta_1(\theta_1-2 \theta_2) - 12 z_1 (6 \theta_1 +5)(6 \theta_1+1) \right) f(z_1, z_2) \\
0 & = & L_2 f(z_1, z_2) \\
\nonumber & = & \left( \theta_2^2-z_2(2\theta_2 - \theta_1+1)(2 \theta_2 - \theta_1) \right) f(z_1, z_2) \ .
\end{eqnarray}
See also \cite{YY} for some related discussion of the Lian-Yau example.

Changing coordinates into the $ (b^2, d) $ affine patch on $ \bb{WP}(2,3,6) $ gives
\[ (z_1, z_2) \mapsto \left( 864 z_1 -1, 144^3 z_1^2 z_2 \right) = (b^2, d) \]
and plugging this into Equations \ref{GDbd1}, \ref{GDbd2} gives
\begin{eqnarray}
0 & = & - \frac{1}{z_2} L_2 F(z_1, z_2) \\
0 & = & \left( - \frac{1}{z_1} L_1 + 1728 L_2 \right) F(z_1, z_2)
\end{eqnarray}
Therefore, as in Section \ref{toriccurves}, the solution spaces of the Picard-Fuchs equations computed via GKZ
and Griffiths-Dwork are identical.

\end{document}